\documentclass{amsart}
\usepackage{graphicx, amssymb, amsthm, amsfonts, latexsym, epsfig, amsmath}
\usepackage{verbatim}
\usepackage{enumerate}
\usepackage{pstricks,pstricks-add,pst-math}

\newtheorem{thm}{Theorem}[section]
\newtheorem{lem}[thm]{Lemma}

\newtheorem{prop}[thm]{Proposition}
\newtheorem{col}[thm]{Corollary}
\theoremstyle{definition}
\newtheorem{defn}[thm]{Definition}
\newtheorem{eg}[thm]{Example}
\newtheorem{rem}[thm]{Remark}

\newif\ifproof
\prooftrue

\DeclareMathOperator{\Int}{int}

\newcommand{\set}[1]{\left\{#1\right\}}

\newcommand{\B}[1]{\ensuremath{\mathbb{#1}}}
\newcommand{\N}{\B{N}}

\newcommand{\ra}{\rightarrow}

\renewcommand{\include}{\input}

\newcommand{\w}{\omega}

\renewcommand{\d}{\delta}

\newcommand{\nowt}{\emptyset}

\newcommand{\eps}{\varepsilon}

\newcommand{\Z}{\B Z}

\newcommand{\nat}{\mathbb{N}}
\newcommand{\integ}{\mathbb{Z}}

\newcommand{\andysq}{\begin{flushright}$\square$\end{flushright}}
\newcommand{\cl}[1]{\overline{#1}}
\newcommand{\restr}[2]{#1\hspace{-0.15cm}\restriction_{#2}}

\numberwithin{equation}{section}
\begin{document}

\title{Shadowing and Expansivity in Sub-Spaces}

\author[A. D. Barwell]{Andrew D. Barwell}
\address[A. D. Barwell]{Heilbronn Institute for Mathematical Research, University of Bristol, Howard House, Queens Avenue, Bristol, BS8 1SN, UK}
\email[A. D. Barwell]{A.Barwell@bristol.ac.uk}
\author[C. Good]{Chris Good}
\address[C. Good]{School of Mathematics, University of Birmingham, Birmingham, B15 2TT, UK}
\email[C. Good]{cg@for.mat.bham.ac.uk}
\author[P. Oprocha]{Piotr Oprocha}
\address[P. Oprocha]{AGH University of Science and Technology, Faculty of Applied Mathematics, al. A. Mickiewicza 30, 30-059 Krak\'ow, Poland, -- and -- Institute of Mathematics, Polish Academy of Sciences, ul. \'Sniadeckich 8, 00-956 Warszawa, Poland}
\email{oprocha@agh.edu.pl}


\begin{abstract}
We address various notions of shadowing and expansivity for continuous maps restricted to a proper subset of their domain. We prove new equivalences of shadowing and expansive properties, we demonstrate under what conditions certain expanding maps have shadowing, and generalize some known results in this area. We also investigate the impact of our theory on maps of the interval.
\end{abstract}

\maketitle

\section{Introduction}\label{intro}
Pseudo-orbit tracing, or \textit{shadowing}, relates to the stability of orbits in a dynamical system under small perturbations, and has been studied by many authors in a variety of contexts. It has been studied in the context of numerical analysis \cite{Corless, Corless2, Pearson}, at times being cited as a prerequisite to achieving accurate mathematical models, and also as a property in its own right \cite{Chen, Coven, KO2, sakai, Nusse, Pennings, Pil, Walters}. Bowen was one of the first to consider this property in \cite{Bowen}, where he used it in the study of $\omega$-limit sets of Axiom A diffeomorphisms. In \cite{BDG}, we use shadowing to characterize $\omega$-limit sets of tent maps, and, following on from \cite{Barwell}, in \cite{BGOR} we use various forms of shadowing to characterize $\omega$-limit sets of topologically hyperbolic systems. Of particular interest is a property called \textit{h-shadowing}, which we prove is equivalent to shadowing in certain expansive systems (such as shifts of finite type), but is in general a stronger property, and one which allows us to prove when internally chain transitive sets are necessarily $\omega$-limit sets \cite{BGOR}.
The notion of an expanding (or expansive) map is closely related to various dynamical properties. In  \cite{Blokh}, Blokh et al use one notion to characterize $\omega$-limit sets of interval maps; in \cite{PEsakai}, Sakai explores various connections between expansivity and shadowing; and in \cite{Urbanski}, Przytycki and Urba{\'n}ski prove shadowing exists for open maps which are expanding on the whole space.

In this paper, we investigate the relationships between various notions of expansivity, the relationships between various notions of shadowing and examine the role that expansivity plays in shadowing. In Section 2, we analyze notions of expansivity that appear both explicitly and implicitly in the literature. In particular, we show that the natural notion of distance expanding (Definition \ref{uniexp}) and the natural notion of ball expanding (Definition \ref{defn:expandingballs}) are the same (Theorem \ref{thm:exp_equivs}) away from what might be termed \emph{critical points}, i.e. points where the map is not locally open or not locally one-to-one.  We do similar analysis of notions of shadowing in Section 3.
By considering these properties on proper subsets of the space, we are able to extend earlier results from other papers, identifying some subtle changes to the theory in this case. In Section 4, we show how maps with different types of expansivity on subsets of their domain have the various types of shadowing on these sets. This allows us to extend Przytycki and Urbanski's result \cite{Urbanski} that an open expanding map has shadowing in Theorem \ref{CINEshadow} by weakening the assumptions and strengthening the conclusion.  These results are summarized in Figure \ref{familytree}. In Section 5 we focus our attention on interval maps. There are many interval maps which do not have shadowing on the whole interval. Nevertheless, it is often possible to identify regions where shadowing does occur, especially away from the critical points.
Coven et al \cite{Coven} show that tent maps have shadowing precisely when the critical point obeys certain parity rules with respect to its orbit; Chen \cite{Chen} generalizes this, proving that maps conjugate to piecewise linear interval maps with gradient modulus greater than $1$ at non-critical points have shadowing if and only if their critical points obey the \textit{linking property}. Theorem \ref{thm:linear_shad} shows that strong shadowing properties hold in piecewise linear maps with gradient modulus greater than 1 away from the critical points.
Piecewise linear maps with gradient modulus greater than $1$ are expanding in a very strong way on open sets separated from the set of critical points. There are however smooth interval maps which also have shadowing properties (such as the logistic map; see Example \ref{eg:logistic_shad}) but which do not share the nice expanding properties of linear maps. Theorem \ref{misiurewicz type} addresses this.
We end with a few remarks showing that h-shadowing homeomorphisms are only possible on totally disconnected sets.
Throughout the paper we provide examples which help illustrate our theory and definitions.


\section{Expansivity and Expanding Maps}\label{linexp}

In this section we explore two types of expansion in maps: \textit{expansive} maps, in which distinct points must, for some iterate of the map, be a set distance apart, and \textit{expanding} maps, for which distinct points which start close together are mapped further apart. The properties defined have been studied extensively, and can be found in many texts, including \cite{AH, sakai, Urbanski, PEsakai, Walters}. We aim to demonstrate both differences and similarities between these properties, and to show that there are cases where maps have these properties on proper subsets of their domain only.

For $x\in X$, we say that $f$ is \emph{open at $x$} if for every neighbourhood $U$ of $x$, $f(x)\in\Int(f(U))$, and that $f$ is \emph{locally one-to-one} at $x$ if there is an open set $V\ni x$ such that $\restr{f}{V}$ is injective. For a subset $\Lambda\subseteq X$, we say that $f$ is \emph{open on $\Lambda$} if for every $x\in\Lambda$ $f$ is open at $x$, and $f$ is \emph{locally one-to-one} on $\Lambda$ if for every $x\in \Lambda$ $f$ is locally one-to-one at $x$. Notice that if $f$ is locally one-to-one on $\Lambda$, then it is locally one-to-one on an open neighbourhood of $\Lambda$.

For any open cover $U$ of $\Lambda$, the \emph{Lebesgue number} of this cover is the constant $\delta$ such that for any $x\in \Lambda$, the open $\delta$-neighbourhood around $x$ is contained in some member of the cover.

\begin{rem}\label{rem:openbasis}
Note that $f$ is open on $\Lambda$ if and only if for every $x\in \Lambda$ there is a neighborhood basis $\set{U_i}_{i\geq 0}$ such that $f(U_i)$ is open, for every $i\geq 0$.
\end{rem}

Our (local) definition of openness is consistent with the standard definition of an open map. Namely, by Remark~\ref{rem:openbasis} if
$f$ is open on $X$ then $f(U)$ is open for every open set $U$.

We say that $f$ is \emph{positively expansive} (with expansive constant $b>0$) if for any $x,y\in X$ the condition
\[
d(f^n(x),f^n(y))<b \hspace{0.5cm}\mbox{ for every }n\in\integ,\ n\geq 0
\]
implies that $x=y$. Moreover, if $f$ is a surjective map it is said to be \textit{$c$-expansive} (with expansive constant $b'>0$) if for any $x,y\in X$ and any full orbits $\set{x_m}_{m\in\Z}$ and $\set{y_n}_{n\in\Z}$ through $x$ and $y$ respectively the condition
\[
d(x_n,y_n)<b' \hspace{0.5cm}\mbox{ for every }n\in\integ
\]
implies that $x=y$. A continuous map $f$ is said to be \textit{topologically hyperbolic} if it is $c$-expansive and has the shadowing property.

There is a large class of topologically hyperbolic maps. The classical example is an Axiom~A diffeomorphism restricted to its
non-wandering set (see \cite{Bowen} for example). Another important class are shifts of finite type (one or two-sided).
The reader is referred to \cite{AH} for a more complete exposition on properties of c-expansive and topologically hyperbolic maps (called TA-maps in that text).

Many authors refer to maps which are expanding on the whole space \cite{sakai, Pil, Urbanski, PEsakai}; there are many situations however where a map will be expanding only in a local sense i.e. only on a subset of the space (such as interval maps which are not expanding on any neighbourhood of a critical point). To tackle this issue, we consider expansion on subsets, and obtain results in this context; the general case of expansion on the whole space is always a natural consequence.

\begin{defn}\label{uniexp}
Let $(X,d)$ be a compact metric space, and $f\colon X\rightarrow X$ be continuous. If there are $\delta>0$, $\mu>1$ such that
$d(f(x),f(y))\geq \mu d(x,y)$ provided that $x,y\in \Lambda\subseteq X$ and $d(x,y)<\delta$ then we say that $f$ is \textit{expanding on $\Lambda$.} If there is a neighbourhood $U$ of $\Lambda$ such that the property holds for every $x,y\in U$, we say that $f$ is \textit{expanding on a neighbourhood of $\Lambda$}. In the case that $\Lambda=X$ we simply say that $f$ is expanding.
\end{defn}

\begin{rem} Przytycki and Urba{\'n}ski explore expanding maps in \cite{Urbanski}, referring to the property as \textit{distance expanding}. They also define a property called \textit{expanding at $\Lambda$}, which is equivalent to our notion of expanding on a neighbourhood of $\Lambda$.
If $f$ is expanding on $\Lambda$ then for each $x\in \Lambda$ there is an open set $U\ni x$ such that $\restr{f}{U\cap \Lambda}$ is one-to-one. Furthermore, if $f$ is expanding on the invariant set $\Lambda$ then it is easy to see that $f$ is positively expansive on $\Lambda$. It is also known that a positively expansive map is expanding with respect to some equivalent metric (see \cite{Urbanski}).
\end{rem}

The following property (for the case $\Lambda=X$) was introduced implicitly in \cite{Urbanski} in the proof that expanding open maps have shadowing. In fact, it implies stronger shadowing properties (see Theorem \ref{CINEshadow}). Here $B_\eps(x)$ denotes the open ball $\{y:d(x,y)<\eps\}$ about $x$.

\begin{defn}\label{defn:expandingballs}
For a compact metric space $X$ with metric $d$ and a subset $\Lambda\subseteq X$ we say that a continuous map $f\colon
X\rightarrow X$ is \textit{ball expanding on $\Lambda$} if there are a $\mu>1$ and a $\nu>0$ such that for every $x\in\Lambda$ and every $\varepsilon<\nu$ we have that $B_{\mu\varepsilon}(f(x))\subseteq f\left(B_{\varepsilon}(x)\right)$.
\end{defn}

\begin{thm}\label{thm:exp_equivs}
Let $(X,d)$ be a compact metric space, $f\colon X\rightarrow X$ be continuous and let $\Lambda\subseteq X$ be closed (but not necessarily invariant). The following are equivalent:
 \begin{enumerate}
 \item \label{open_and_exp} $f$ is open on a neighbourhood of $\Lambda$ and expanding on a neighbourhood of $\Lambda$,
 \item \label{top_expanding_and_1-1} $f$ is ball expanding on a neighbourhood of $\Lambda$ and locally one-to-one on (a neighbourhood of) $\Lambda$.
 \end{enumerate}
\end{thm}

\ifproof\begin{proof} \eqref{open_and_exp} $\Longrightarrow$ \eqref{top_expanding_and_1-1}: Let $Q$ be the neighbourhood of $\Lambda$ on which $f$ is open and expanding, and let $\mu>1$, $\delta>0$ be the constants as given in the definition of expanding. By normality, there is some open set $W$ such that $\Lambda\subseteq W\subseteq\cl{W}\subseteq Q$ (if $\Lambda=X$ take $\Lambda = W = \cl{W} = Q$). Clearly $f$ is one-to-one on $Q$ and therefore one-to-one on $\Lambda$.

Fix $x\in \cl{W}$, then there is some $\zeta<\delta$ such that $f$ is open and expanding on $U=B_{\zeta}(x)\subseteq Q$. Furthermore there is an $\eta=\eta(x)>0$ such that $B_{\mu\eta}(f(x))\subset\Int(f(U))$ since $f$ is open at $x$. Fix any $\rho<\eta$. We claim that $B_{\mu \rho}(f(x))\subseteq f(B_{\rho}(x))$. To see this, let $V=U \cap f^{-1}\bigl(B_{\mu \rho}(f(x))\bigr)$ and notice that $f(V)=B_{\mu \rho}(f(x))$. Suppose that $V\nsubseteq B_{\rho}(x)$. Then there is $y\in V\setminus B_{\rho}(x)$, so $d(x,y)\geq \rho$ and $x,y\in U \subseteq Q$, thus $d(f(x),f(y))\geq \mu \rho$ and we get that $f(y)\not\in B_{\mu \rho}(f(x))$ which is impossible, so the claim holds.

Let $U'=f^{-1}\bigl(B_{\mu\eta/2}(f(x))\bigr)\cap U$. Then $x\in U'$ and we can take $\nu=\nu(x)<\eta/2$ so that $B_{\nu}(x)\subseteq U'$.

Take any $z\in B_{\nu}(x)$ and $\eps\leq\nu$ so that $B_{\eps}(z)\subseteq B_{\nu}(x)$, then $f(z)\in f(U')= B_{\mu\eta/2}(f(x))$ and so
\[
B_{\mu\eps}(f(z))\subseteq B_{\mu\eta}(f(x)) \subseteq \Int f(U).
\]
Since $z\in U$, a similar argument shows that $B_{\mu \eps}(f(z))\subseteq f(B_\eps(z))$. In other words, for $x\in\cl{W}$ and $z\in X$ we have that
\begin{align}
B_\eps(z)\subseteq B_{\nu(x)}(x)\ \Longrightarrow\ B_{\mu \eps}(f(z))\subseteq f(B_\eps(z)).\label{ex_balls_z}
\end{align}
Note that $\cl{W}$ is compact and $\nu(x)$ is well defined for every $x\in \cl{W}$, so there are $x_1,\ldots, x_s$ such that
\[
\cl{W} \subseteq \bigcup_{i=1}^s B_{\nu(x_i)/2}(x_i).
\]
Denote $\xi = \min_{i} \nu(x_i)/2$, fix any $\eps<\xi$ and any $x\in \cl{W}$. There is $i$ such that $x\in B_{\nu(x_i)/2}(x_i)$ and so $B_\eps(x)\subseteq B_{\nu(x_i)}(x_i)$. Hence by \eqref{ex_balls_z}, $B_{\mu\eps}(f(x))\subseteq f(B_{\eps}(x))$.

\eqref{top_expanding_and_1-1} $\Longrightarrow$ \eqref{open_and_exp}: Let $W$ be the neighbourhood of $\Lambda$ on which $f$ is ball expanding and let $\mu$, $\nu$ be as given in Definition \ref{defn:expandingballs}. Certainly for every $x\in\Lambda$ there is a $\zeta(x)$ such that $f$ is one-to-one on $B_{\zeta(x)}(x)$, and the collection of such neighbourhoods cover $\Lambda$. Take $\beta$ to be their Lebesgue number and let $\eps:=\min\{\beta,\nu\}$, then $f$ is one-to-one on $B_{\eps}(x)$ for every $x\in\Lambda$.

Now consider a cover of $\Lambda$ consisting of $\eps/3$-neighbourhoods of points in $\Lambda$, take a finite subcover $\{B_{\eps/3}(x_i)\ :\ 1\leq i\leq n\}$, and let $U=W\cap\bigcup_{i\leq n} B_{\eps/3}(x_i)$. Let $\delta=\eps/3$ and fix any $x,y\in U$ with $\eta=d(x,y)<\delta$, then for some $i\leq n$, $d(x,x_i)<\eps/3$.

Suppose that $d(f(x),f(y))<\mu d(x,y)=\mu\eta$. Then we have that
\[f(y)\in B_{\mu\eta}(f(x))\subseteq f(B_{\eta}(x))\]
since $f$ is ball expanding at $x$, and $\eta<\nu$. But $y\notin B_{\eta}(x)$, so there is a $z\in B_{\eta}(x)$ for which $f(z)=f(y)$. Since both $y$ and $z$ are in $B_{\eps}(x_i)$ this contradicts the fact that $f$ is one-to-one on $B_{\eps}(x_i)$. Thus $d(f(x),f(y))\geq\mu d(x,y)$ for every $x,y\in U$ with $d(x,y)<\delta$, and hence $f$ is expanding on $U$.

To see that $f$ is open on $W$, take any $x\in W$ and any $0<\eps<\nu$. Then $B_{\mu\eps}(f(x))\subseteq f(B_{\eps}(x))$, which implies that $f(x)\in\Int f(B_{\eps}(x))$.
\end{proof}\fi

\begin{eg}
  Suppose that $f:[0,1]\to[0,1]$ is a piecewise linear map such that $f(0)=1/2$ and the gradient on $[0,1/8]$ is $3/2$. Then $f$ is neither open at $0$ nor ball expanding on $\Lambda=[0,1/16]$ but it is locally one-to-one and expanding on a neighbourhood of $\Lambda$.
\end{eg}

The distinction between expanding on $\Lambda$ and expanding on a neighbourhood of $\Lambda$ suggests the following intermediate property:
               \begin{enumerate}
               \item[($\star$)] there are $\delta>0$, $\mu>1$ such that $d(f(x),f(y))\geq \mu d(x,y)$ provided that $x\in\Lambda$ and $d(x,y)<\delta$.
               \end{enumerate}
The property ($\star$) is intermediate between expanding on $\Lambda$ and on a neighbourhood, namely it is immediate that maps which are expanding on a neighbourhood of $\Lambda$ have ($\star$) and that ($\star$) implies expanding on $\Lambda$. The proof of the following is similar to that of Theorem \ref{thm:exp_equivs} and is left to the reader.

\begin{thm}\label{expequiv1}
Let $(X,d)$ be a compact metric space, $f\colon X\rightarrow X$ be continuous and $\Lambda\subseteq X$ be closed.
\begin{enumerate}
\item If $f$ is  open on $\Lambda$ and has ($\star$), then $f$ is ball expanding on $\Lambda$.
\item If $f$ is ball expanding and locally one-to-one on $\Lambda$, then $f$ has ($\star$).
\end{enumerate}
\end{thm}

The next two examples illustrate the differences between maps which are expanding and those which are ball expanding.

\begin{eg}\label{eg:Cantoreg}
There is a continuous function from the Cantor set to itself that is expanding but not ball expanding. Our Cantor set $X$ is the subset of $[-1,1]$ consisting of the union of the middle third Cantor set on $[0,1]$ and its left shift by $-1$, with the usual metric. For each $0<n$, let  $C_n=X\cap[2/3^{n}, 1/3^{n-1}]$ and $C_{-n}=X\cap[-1/3^{n-1},-2/3^{n}]$, so that $X=\{0\}\cup\bigcup_{n\in\Z\setminus \set{0}}C_n$.
For $n>0$, let $C_n^+=C_n\cap [8/3^{n+1},1/3^{n-1}]$ be the right hand half of $C_n$ and $C_n^-=[2/3^{n},7/3^{n+1}]$ be the left hand half, so that $C_n=C_n^-\cup C_n^+$. Clearly $3C_n=C_{n-1}$ for $1<|n|$ and $C_n^+$ and $C_n^-$ are both isometric copies, indeed translations, of $C_{n+1}$. Our function $f$ fixes $0$ and expands each $C_n$ by a multiple of $3$ and then translates it so that it embeds it into $X$ in the following way. $C_{-3}$ and $C_3$ are expanded by a factor of 9, all other $C_n$s are expanded by a factor of 3 and
$f(C_{-1})=X\cap[-1,0]$, $f(C_1)=X\cap [0,1]$,
$$
f(C_n)=
\begin{cases}
  C_1&\text{ if }n=\pm2,\pm3\\
  C_{n-2}^+&\text{ if }n>3,\\
  C_{|n|-2}^-&\text{ if }n<-3.\\
\end{cases}
$$
Clearly $f$ is expanding, with $\delta=1/9$ and $\mu=3$, for example, in Definition \ref{uniexp}. However, for $n\geq4$, $f\big(X\cap(-2/3^n,2/3^n)\big)=X\cap [0,1/3^{n-3}]$,  so that for any $\eps<1/27$, if $\delta>0$, then $B_\delta(f(0))$ is not a subset of $f(B_\eps(0))$. Hence there is no $\mu>1$ such that $B_{\mu\eps}\big(f(0)\big)\subseteq f\big(B_\eps(0)\big)$, i.e. $f$ is not ball expanding.\andysq
\end{eg}

\begin{eg}\label{eg:equiv_tent}
The full tent map $T_2$ is an example of a map which is ball expanding on $[0,1]$ but not expanding and not locally one-to-one.

Indeed, since $T_2$ is unimodal and the image of the critical point and the end points is an endpoint, $T_2$ is open on $[0,1]$.  Since the gradient has modulus $2$ except at 1/2, $T_2$ is therefore easily seen to be ball expanding. However, as it is open, but not one-to-one on any neighbourhood of the critical point, $T_2$ is not expanding on $[0,1]$.
\andysq
\end{eg}


\section{Shadowing on Subspaces}\label{shad}

Let $X$ be a compact metric space and $f\colon X\rightarrow X$ continuous. For $\eps>0$, the (finite or infinite) sequence $\{x_0,x_1,\ldots\}\subseteq X$ is an  $\eps$-\emph{pseudo-orbit} if $d(f(x_n),x_{n+1})<\eps$, for all $n\ge 0$. The sequence (when infinite) is an \emph{asymptotic pseudo-orbit} if $d(f(x_n),x_{n+1}) \to 0$ as $n\to \infty$ and an \emph{asymptotic} $\eps$-\emph{pseudo-orbit} if both conditions hold.

Let $\eps>0$, and let $K$ be either $\nat$ or $\{0,1,\ldots,k-1\}$ for some $k\in\nat$. The sequence $\{y_n\}_{n\in K}$ $\eps$-\emph{shadows} the sequence $\{x_n\}_{n\in K}$ if and only if for every $n\in K$, $d(y_n,x_n)<\eps$. Furthermore, we say that the sequence $\{y_n\}_{n\in\nat}$ \emph{asymptotically shadows} the sequence $\{x_n\}_{n\in\nat}$ if and only if $\lim_{n\rightarrow\infty}d(x_n,y_n)=0$. If both conditions hold simultaneously, we say that $\{y_n\}_{n\in\nat}$ \emph{asymptotically $\eps$-shadows} the sequence $\{x_n\}_{n\in\nat}$. If $y_n=f^n(y)$ for every $n\in\nat$ then we say that the point $y$ shadows (in whichever sense is appropriate) the sequence $\{x_n\}_{n\in\N}$.

To compliment our treatment of expansive properties, we introduce many of the following definitions with respect to the given set, as well as in a general form.

The standard version of shadowing is the following, which appeared in \cite{Bowen}, where it was used in the study of $\omega$-limit sets of Axiom A diffeomorphisms.

\begin{defn}
Let $(X,d)$ be a compact metric space, $f\colon X\rightarrow X$ be continuous and let $Y$ be a subset of $X$. We say that $f$ has the \textit{pseudo-orbit tracing property on $Y$} (or \textit{shadowing on $Y$}) if and only if for every $\varepsilon>0$ there is $\delta>0$ such that every infinite $\delta$-pseudo-orbit in $Y$ is $\eps$-shadowed by a point $y\in X$.

If this property holds on $Y = X$, we simply say that $f$ has \emph{shadowing}.
\end{defn}

\begin{rem}\label{lem:finite_shadowing}
It is easy to see that $f$ has shadowing if and only if for every $\eps>0$ there is a $\delta>0$ such that every finite $\delta$-pseudo-orbit is $\eps$-shadowed.
\end{rem}


\begin{defn} Let $(X,d)$ be a compact metric space, $f\colon X\rightarrow X$ be continuous and let $Y$ be a subset of $X$. We say that $f$ has
\emph{limit shadowing on $Y$} if and only if for any asymptotic pseudo-orbit $\set{x_n}_{n \in \N} \subseteq Y$ there is a point $y\in X$ which asymptotically shadows $\set{x_n}_{n \in \N}$.

If this property holds on $Y = X$, then we say that $f$ has \emph{limit shadowing}.
\end{defn}

The definition of limit shadowing was extended in \cite{sakai} to a property called \textit{s-limit shadowing} (Definition \ref{defn:s-lim-shad}), to account for the fact that many systems have limit shadowing but not shadowing \cite{KO2,Pil}.

\begin{defn}\label{defn:s-lim-shad}
Let $(X,d)$ be a compact metric space, and $f\colon X\rightarrow X$ be continuous.
We say that \emph{$f$ has s-limit shadowing on $Y \subseteq X$} if and only if
for every $\eps > 0$ there is $\delta > 0$ such that the following two
conditions hold:
\begin{enumerate}
    \item for every $\delta$-pseudo-orbit $\set{x_n}_{n\in \N}\subseteq Y$ of $f$, there is $y \in X$ such that $y$ $\eps$-shadows $\{x_n\}_{n\in\N}$, and
    \item for every asymptotic $\delta$-pseudo-orbit $\set{z_n}_{n\in \N}\subseteq Y$ of $f$, there is $y \in X$ such that $y$ asymptotically $\eps$-shadows $\{z_n\}_{n\in\N}$.
\end{enumerate}
In the special case $Y=X$ we say that $f$ has \emph{s-limit shadowing}.
\end{defn}

\begin{eg}
  Let $X=[0,1]\cup\{-1/2^n:n\ge1\}$. Let $f\colon X\to X$ be any homeomorphism such that $f(x)=x$ for $x=1$ or $x\le0$ and $f(x)<x$ for $x\in (0,1)$.
  We claim that $f$ has shadowing and limit shadowing, but that it does not have s-limit shadowing.

  Note first that if $\delta<1/2^n$, for some $1\le n$,  then any $\delta$-pseudo orbit, $\{x_i\}_{i\ge0}$, is either a constant sequence or is contained in $[-1/2^n,1]$. In the second case, there is a $\delta$-pseudo
orbit $\{y_i\}_{i\ge0}\subseteq[0,1]$ such that $d(x_i,y_i)\leq 1/2^n$. If $\{x_i\}_{i\ge0}$ is an asymptotic pseudo-orbit in $X$, then either it is eventually constant, or there is an asymptotic pseudo orbit $\{y_i\}_{i\ge0}\subseteq[0,1]$ such that $d(x_i,y_i)\to0$ as $i\to\infty$. But it is also well-known that the restriction of $f$ to $[0,1]$ has shadowing and limit shadowing. It follows that $f$ does as well.

To see that $f$ does not have s-limit shadowing, let $\eps=1/4$ and choose any $\delta>0$. Fix $N>0$ such that $f^N(1/2)<\delta$ and $2^{-N}<\delta$. Then the sequence
$$
\{x_i\}_{i\ge0}=\left\{\frac{1}{2}, f\left(\frac{1}{2}\right), \ldots , f^N\left(\frac{1}{2}\right),0, -\frac{1}{2^N}, -\frac{1}{2^N}, \ldots\right\}
$$
is both a $\delta$-pseudo orbit and an asymptotic pseudo-orbit. But now, if $z\in X$ satisfies $\lim_{i\to\infty}d(f^i(z),x_i)=0$ then $z=-1/2^N$ and so it does not $\eps$-trace
our $\delta$-pseudo orbit $\{x_i\}_{i\ge0}$. Indeed, $f$ does not have s-limit shadowing.

We conjecture that one might extend this example to an interval map but the details of the proof seem convoluted.\andysq
\end{eg}

Walters \cite{Walters} showed that a shift space is of finite type if and only if it has shadowing. The following definition was introduced in \cite{BGOR} and is motivated by the fact that shifts of finite type actually possess a stronger shadowing property, which happens to coincide with shadowing in shift spaces (but not in other systems -- see Example \ref{weakmixeg}). Later, we will show that h-shadowing is satisfied by various interval maps on regions excluding local extrema.

\begin{defn}
Let $(X,d)$ be a compact metric space, and $f\colon X\rightarrow X$ be continuous. We say that $f$ has \textit{h-shadowing on $Y\subseteq X$} if and only if for every $\varepsilon>0$ there is a $\delta>0$ such that for every finite $\delta$-pseudo-orbit $\{x_0,x_1,\ldots,x_m\}\subseteq Y$ there is $y\in X$ such that $d(f^i(y),x_i)<\varepsilon$ for every $i<m$ and $f^m(y)=x_m$.

If $Y=X$ then we simply say that $f$ has \emph{h-shadowing}.
\end{defn}

It is easy to see that every map with h-shadowing has shadowing. The converse is not true however, as Example \ref{weakmixeg} shows.
The following theorem relates h-shadowing, s-limit shadowing and limit shadowing. (In \cite{BGOR} it was shown that, if $\Lambda\subseteq f(\Lambda)\subseteq X$, in particular if $f$ is onto, and $f$ has s-limit shadowing on $\Lambda$ then $f$ has also limit shadowing on $\Lambda$.)

\begin{thm}\label{h-shad_and_s-limit_shad}
Let $(X,d)$ be a compact metric space, $f\colon X\rightarrow X$ be continuous and suppose that $\Lambda\subseteq X$ is closed.
\begin{enumerate}
\item If $f$ has s-limit shadowing on $\Lambda$, then $f$ has limit shadowing on $\Lambda$.\label{h-shad_and_s-limit_shad0}
    \item If there is an open set $U$ such that $\Lambda\subseteq U$ and $f$ has h-shadowing on $U$, then $f$ has s-limit shadowing on $\Lambda$. \label{h-shad_and_s-limit_shad1}
    \item If $\Lambda$ is invariant and $\restr{f}{\Lambda}$ has h-shadowing then $\restr{f}{\Lambda}$ has s-limit shadowing and limit shadowing. \label{h-shad_and_s-limit_shad2}
    \item If $f$ has h-shadowing then $f$ has s-limit shadowing and limit shadowing. \label{h-shad_and_s-limit_shad3}
\end{enumerate}
\end{thm}
\ifproof\begin{proof}
\eqref{h-shad_and_s-limit_shad0}:
To prove limit shadowing, take any asymptotic pseudo-orbit $\set{z_n}_{n\in \N}$ in $\Lambda$. Fix $\eps>0$ and let $\delta$ be provided for $\eps$
by s-limit shadowing, and let $\gamma<\delta$ be provided for $\delta$ by the same condition. There is $K$ such that $\set{z_n}_{n\geq K}$ is a $\gamma$-pseudo-orbit, so it is asymptotically $\delta$-shadowed by a point $z$. It follows that the $w$-limit set, $\w(z)$, is a subset of $\Lambda$. By \cite[Theorem 8.7]{FurBook}, there exists a minimal subset $M$ of $\w(z)$ and a point $y\in M$ such that $\liminf_{j\to\infty}d(f^j(z),f^j(y))=0$. But $f|_M$ is onto, so there exists a point $x$ such that $f^K(x)=y$. There is also $N>0$
such that $d(f^N(y),f^N(z))<\delta/2$ and $d(f^N(z),z_{N+K})<\delta/2$. Therefore the sequence
$$
\xi=\set{x,f(x),\ldots, f^{K+N-1}(x), z_{N+K}, z_{N+K+1},\ldots}
$$
is an asymptotic $\delta$-pseudo-orbit in $\Lambda$. Now, it is enough to use s-limit shadowing obtaining a point which asymptotically $\eps$-traces $\xi$, and as a result asymptotically traces $\set{z_n}_{n\in \N}$.

\eqref{h-shad_and_s-limit_shad1}: Since every map with h-shadowing has shadowing, the first half of the definition of s-limit shadowing is satisfied trivially.

So fix $\eps>0$ such that $B(\Lambda,3\eps)\subseteq U$ and denote $\eps_n=2^{-n-1}\eps$. By the definition of $h$-shadowing there are
$\set{\delta_n}_{n\in \N}$ such that every $\delta_n$-pseudo-orbit in $U$ is $\eps_n$-traced (with exact hit at the end). Fix any
infinite $\delta_0$-pseudo-orbit $\set{x_n}_{n\in\N}\subseteq \Lambda$ such that $\lim_{n\ra \infty} d(f(x_n),x_{n+1})=0$. There is an
increasing sequence $\set{k_i}_{i\in \N}$ such that $\set{x_n}_{n=k_i}^\infty$ is an infinite $\delta_i$-pseudo-orbit and obviously $k_0=0$. Note that if $w$ is a point such that  $f^{k_i}(w)=x_{k_i}$ then the sequence
\[
w,f(w),\ldots, f^{k_i}(w), x_{k_i+1},\ldots,x_{k_{i+1}}
\]
is a $\delta_{i}$-pseudo-orbit.

Let $z_0$ be a point which $\eps_0$-shadows the $\delta_0$-pseudo-orbit $x_0,\ldots, x_{k_1}$ with exact hit (i.e. such that
$f^{k_1}(z_0)=x_{k_1}$). Notice that $f^j(z_0)\in U$ for $0\leq j\leq k_1$.

For $i\in\nat$, assume that $z_i$ is a point which $\eps_i$-shadows the $\delta_i$-pseudo-orbit
\[
z_{i-1},f(z_{i-1}),\ldots,f^{k_i}(z_{i-1}),x_{k_i+1},\ldots,x_{k_{i+1}}\subseteq U
\]
with exact hit. Then by h-shadowing there is a point $z_{i+1}$ which $\varepsilon_{i+1}$-shadows the $\delta_{i+1}$-pseudo-orbit
\[
z_i,f(z_i),\ldots,f^{k_{i+1}}(z_i),x_{k_{i+1}+1},\ldots,x_{k_{i+2}}\subseteq U
\]
with exact hit. Thus we can produce a sequence $\set{z_i}_{i=0}^\infty$ with the following properties:
\begin{enumerate}[(a)]
\item $d(f^j(z_{i-1}),f^j(z_i)) < \eps_i$ for $j\leq k_{i}$ and $i\geq 1$,\label{c1:hhs}
\item $d(f^j(z_i),x_j) < \eps_i$ for {$k_i< j\leq k_{i+1}$} and $i\geq 0$,\label{c2:hhs}
\item $f^{k_{i+1}}(z_i)=x_{k_{i+1}}$ for $i\geq0$,\label{c3:hhs}
\item $d(f^j(z_i),\Lambda)<\eps$ for $j\leq k_{i+1}$,\label{c4:hhs}
\end{enumerate}
There is an increasing sequence $\set{s_i}_{i\in \N}$ such that the limit $z=\lim_{i\ra \infty} z_{s_i}$ exists.

For any $j,n\in\nat$ there exist $i_0\geq0$ and $m\geq i_0$ such that $k_{i_0}<j\leq k_{i_0+1}$ and $d(f^j(z),f^j(z_{s_m}))<\eps_{n+1}$. So we get
\begin{eqnarray*}
d(f^j(z),x_j)&\leq & d(f^j(z),f^j(z_{s_m}))+d(f^j(z_{i_0}),x_j)+\sum\limits_{i=i_0}^{s_m-1}d(f^j(z_i),f^j(z_{i+1}))\\
&\leq & \eps_{n+1}+\eps_{i_0}+\sum_{i=i_0}^{s_m-1}\eps_{i+1}\\
&\leq & \eps 2^{-n-2}+\sum_{i=i_0}^\infty 2^{-i-1}\eps \; \leq\; \eps (2^{-n-2}+2^{-i_0})\\
&\leq & \eps (2^{-n-2}+1).
\end{eqnarray*}
But we can fix $n$ to be arbitrarily large in that case, which immediately implies that
\begin{eqnarray*}
d(f^j(z),x_j)&\leq & \eps.
\end{eqnarray*}
Furthermore, for any $n$, let $j>k_{n+2}$. There is $i_1\geq n+2$ such that $k_{i_1}<j\leq k_{i_1+1}$ and there is $m>i_1$ such that $d(f^j(z),f^j(z_{s_m}))<\eps_{n+1}$. Then as before we obtain
\begin{eqnarray*}
d(f^j(z),x_j)&\leq & \eps (2^{-n-2}+2^{-i_1})\\
&\leq & \eps (2^{-n-2}+2^{-n-2})=\eps_n.
\end{eqnarray*}
This immediately implies that $\limsup_{j\ra\infty}d(f^j(z),x_j)\leq \eps_n$ which, since $n$ was arbitrary, finally gives $\lim_{j\ra \infty} d(f^j(z),x_j)=0$. This shows that $f$ has s-limit shadowing on $\Lambda$.

\eqref{h-shad_and_s-limit_shad2} and \eqref{h-shad_and_s-limit_shad3} follow directly from \eqref{h-shad_and_s-limit_shad0} and \eqref{h-shad_and_s-limit_shad1} (since $U=\Lambda$ is open in $\Lambda$). \end{proof}\fi

We finish this section by proving a result which shows that provided we can find some iterate of a map which has h-shadowing, we can deduce that the map itself has h-shadowing.
We need the following result, which follows easily from the definitions.

\begin{lem} \label{finite shadowing}
Let $(X,d)$ be a compact metric space, and $f\colon X\rightarrow X$ be continuous. Let $\epsilon>0$ and $n\in\nat$. There is a $\d=\d(n,\eps)>0$ such that if $\{x_0,\ldots,x_n\} $ is a $\d$-pseudo-orbit and $y\in X$ is such that $d(y,x_0)<\d$ then $d(f^k(y),x_k)<\epsilon$ for $k=1,\ldots,n$.
\end{lem}

\begin{thm}\label{iterateshadow}
Let $(X,d)$ be a compact metric space, and $f\colon X\rightarrow X$ be continuous. If $\Lambda$ is a closed set such that $f(\Lambda)\supseteq \Lambda$ then the following conditions are equivalent:
\begin{enumerate}
\item $f$ has $h$-shadowing on $\Lambda$,\label{iterateshadow:1}
\item $f^n$ has $h$-shadowing on $\Lambda$ for some $n\in\nat$,\label{iterateshadow:2}
\item $f^n$ has $h$-shadowing on $\Lambda$ for all $n\in\nat$,\label{iterateshadow:3}
\end{enumerate}
\end{thm}
\ifproof\begin{proof}
Implication from \eqref{iterateshadow:3} to \eqref{iterateshadow:2} is trivial.

Implication from \eqref{iterateshadow:1} to \eqref{iterateshadow:3} is also obvious, since for any $\delta>0$ and $n>0$ if $\set{y_0,y_1,\ldots, y_m}$ is $\delta$-pseudo-orbit for $f^n$ then the sequence
\[
y_0, f(y_0),\ldots, f^{n-1}(y_0), y_1, f(y_1), \ldots, f^{n-1}(y_{m-1}), y_m
\]
is $\delta$-pseudo-orbit for $f$.

For the proof of the last implication fix $\epsilon>0$ and suppose that $f^n$ has $h$-shadowing on $\Lambda$ for some $n\in\nat$. By Lemma~\ref{finite shadowing} there is an $\epsilon'>0$ such that if $\{x_0,\ldots,x_n\}\subseteq \Lambda$ is an $\epsilon'$-pseudo-orbit and $y\in X$ is such that $d(y,x_0)<\epsilon'$ then $d(f^k(y),x_k)<\epsilon$ for $k=1,\ldots,n$.

By $h$-shadowing there is a $\delta>0$ such that every $\delta$-pseudo-orbit of $f^n$ is $\epsilon'$-shadowed by a point in $X$ which hits the last element of the pseudo-orbit. Again by Lemma~\ref{finite shadowing} (with $y=x_0$), there is a $\gamma<\frac{\delta}{n}$ such that whenever $\{x_0,\ldots,x_n\}$ is a $\gamma$-pseudo-orbit for $f$ we have that $d(f^i(x_0),x_i)<\delta$ for $i=1,\ldots,n$.

Let $\{x_0,\ldots,x_m\}\subseteq \Lambda$ be any $\gamma$-pseudo-orbit for $f$, and write $m=jn+r$ for some $j\geq0$ and some $r<n$. Since $f$ is surjective on $\Lambda$ (i.e. $\Lambda \subseteq f(\Lambda)$) there is a point $z\in \Lambda$ such that $f^{n-r}(z)=x_0$. Then $\{z,f(z),\ldots,f^{n-r}(z),x_1,\ldots,x_m\}\subseteq \Lambda$ is a $\gamma$-pseudo-orbit for $f$, which we enumerate obtaining the sequence $\{y_0,\ldots,y_{(j+1)n}\}$. We now claim that $\{y_0,y_n,y_{2n},\ldots,y_{(j+1)n}\}$ is a $\delta$-pseudo-orbit for $f^n$. Indeed, $\{y_0,\ldots,f^{n-r}(y_0)=y_{n-r},\ldots,y_n\}$ is a $\gamma$-pseudo-orbit (of length $n+1$) for $f$ and so $d(f^n(y_0),y_n)<\delta$. Similarly we have $d(f^n(y_{kn}),y_{(k+1)n})<\delta$ for $1\leq k\leq j$.

By $h$-shadowing of $f^n$ there is $u$ such that $d(f^{kn}(u),y_{kn})<\epsilon'$ for $k=0,1,\ldots,j+1$ and $f^{(j+1)n}(u)=y_{(j+1)n}$. Thus by the definition of $\epsilon'$ we have that $d(f^{kn+i}(u),y_{kn+i})<\epsilon$ for $k=0,\ldots,j+1$ and for $i=0,\ldots,n-1$. So the point $u$ $\epsilon$-shadows the $\gamma$-pseudo-orbit $\{y_0,\ldots,y_{(j+1)n}\}=\{z,f(z),\ldots,f^{n-r}(z)=x_0,x_1,\ldots,x_m\}$, and consequently the point $w=f^{n-r}(u)$ $\epsilon$-shadows the $\gamma$-pseudo-orbit $\{x_0,\ldots,x_m\}$, with $f^m(w)=f^{(j+1)n}(u)=y_{(j+1)n}=x_m$.
\end{proof}\fi


\section{Expansivity and Shadowing}\label{sec:exp_and_shad}

In this section we show that the various shadowing properties we have discussed are seen (and many are indeed equivalent) for certain expanding maps.

The following result appears in \cite{BGOR}, and we omit the proof here.

\begin{prop}\label{thm:LimS}
Let $(X,d)$ be a compact metric space and let $f\colon X\rightarrow X$ be continuous.
\begin{enumerate}
\item If $f$ is positively expansive then $f$ has shadowing if and only if $f$ has h-shadowing;
\item If $f$ is $c$-expansive then $f$ has shadowing if and only if $f$ has s-limit shadowing.
\end{enumerate}
\end{prop}

Corollary \ref{hshad_and_limshad} follows immediately from Theorem \ref{h-shad_and_s-limit_shad} and Proposition \ref{thm:LimS}.

\begin{col}\label{hshad_and_limshad}
Let $(X,d)$ be a compact metric space, $f\colon X\rightarrow X$ be continuous and positively expansive.
Then $f$ has shadowing if and only if it has h-shadowing if and only if it is s-limit shadowing. In this case $f$ also has limit shadowing.
\end{col}

The following theorem extends Przytycki and Urbanski's result that open, expanding maps have shadowing \cite[Corollary 3.2.4]{Urbanski}.

\begin{thm}\label{CINEshadow}
Let $(X,d)$ be a compact metric space, $f\colon X\rightarrow X$ be continuous, and let $M\subseteq X$. If $f$ is ball expanding on $M$ then $f$ has h-shadowing on $M$.
\end{thm}
\ifproof\begin{proof}
Let $\varepsilon>0$, let $\mu,\nu$ be as given in Definition \ref{defn:expandingballs}, let $\varepsilon'=\min\{\varepsilon,\nu\}$ and let $\delta=(\mu-1)\varepsilon'$. Then for every $x\in M$
\begin{align}
B_{\varepsilon'+\delta}(f(x))\subseteq B_{\mu\varepsilon'}(f(x)) &\subseteq f\left(B_{\varepsilon'}(x)\right)\label{hash}
\end{align}
Suppose that $\{x_0,\ldots,x_m\}\subseteq M$ is a $\delta$-pseudo-orbit. Notice that by (\ref{hash}) we have $B_{\varepsilon'+\delta}(f(x_i))\subseteq
f\left(B_{\varepsilon'}(x_i)\right)$ for $i=0,1,\ldots,m-1$, so
\begin{equation}
B_{\varepsilon'}(x_{i+1})\subseteq f\left(B_{\varepsilon'}(x_i)\right)\quad\quad \text{ for } i=0,1,\ldots,m-1.\label{con:ball_inclusions}
\end{equation}
Let $J_0=B_{\varepsilon'}(x_0)$ and then define inductively $J_i=f^{-i}\bigl(B_{\varepsilon'}(x_i)\bigr)\cap J_{i-1}$.

Clearly the $J_i$ are nested, and by \eqref{con:ball_inclusions} we can prove by induction that $f^i(J_i)=B_{\varepsilon'}(x_i)$, since
$$
B_{\eps'}(x_i) \supseteq f^i(J_i)\supseteq f^i(J_{i-1})\supseteq f(B_{\eps'}(x_{i-1}))
\supseteq B_{\eps'}(x_i).
$$

In particular, $f^i(J_m)\subseteq B_{\varepsilon'}(x_i)$, for $i=0,1,\ldots,m$ and $f^m(J_m)=B_{\eps'}(x_m)$, thus there is a point $y\in J_m$ such that $f^i(y)\in B_{\varepsilon'}(x_i)$ and for which $f^m(y)=x_m$.
\end{proof}\fi

The following corollary is now immediate from Theorems \ref{thm:exp_equivs} and \ref{CINEshadow}.

\begin{col}\label{col:uniexp_open_hshad}
Let $(X,d)$ be a compact metric space, and $f\colon X\rightarrow X$ be continuous.
\begin{enumerate}
 \item If $f$ is ball expanding, then $f$ has h-shadowing.
  \item If $f$ is open and expanding, then $f$ has h-shadowing.
\end{enumerate}
\end{col}

\begin{center}
\begin{figure}[ht]
    \begin{center}
        \includegraphics[width=0.99\textwidth]{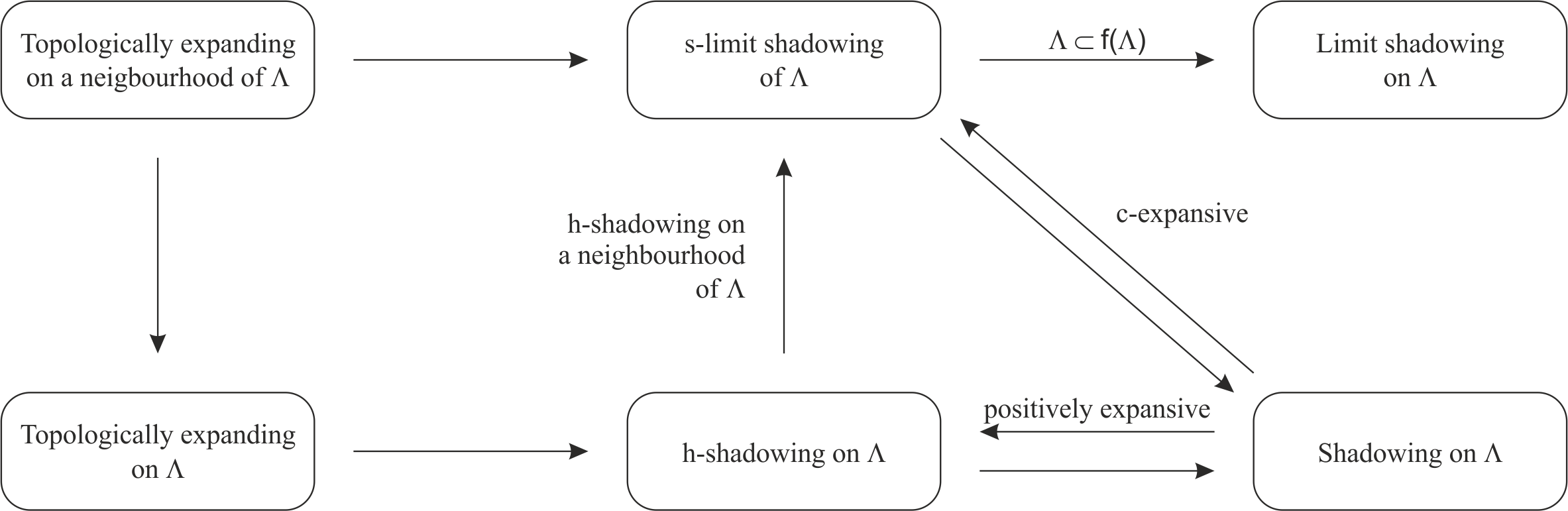}
    \end{center}
    \caption{\small{The relationships between expansivity and shadowing properties of a continuous map $f$ on a subset $\Lambda$ of a compact metric space.}}\label{familytree}
\end{figure}
\end{center}

Figure~\ref{familytree}
summarizes the situation for continuous functions of a compact metric space.

\section{Expansivity and Shadowing in Interval Maps}

We now consider our results in the context of interval maps. Let $f\colon [0,1]\to[0,1]$ be continuous. By a  \emph{critical point} of $f$ we mean a point at which the map fails to be locally one-to-one. We denote the set of critical points by ${\mathcal C}$.

\begin{rem}\label{expintmap} Let $f\colon [0,1]\rightarrow [0,1]$ be continuous. Notice that $\mathcal C$ is a closed subset of $[0,1]$. If $U$ is an open set that contains a critical point, then $f$ is not expanding on $U$, since expanding maps are one-to-one.
On the other hand, if $U$ is open and disjoint from $\mathcal C$, then $U$ can be written as a countable union of disjoint {open connected subsets of $[0,1]$} on which $f$ is one-to-one {(these subsets are open intervals, except maybe at most two of them containing endpoints). Therefore, if $U$ does not} contain either of the end points $0$ or $1$, $f$ is open on $U$. It follows that if $f$ is expanding on an open set $U$ that does not contain any critical points or end points, then $f$ is open on $U$. So if $\Lambda$ is a closed subset of $[0,1]$ that is disjoint from $\mathcal C\cup\{0,1\}$ and $f$ is expanding on a neighbourhood of $\Lambda$, then $f$ is open on a neighbourhood $U$ of $\Lambda$, and by Theorem \ref{thm:exp_equivs} we have that $f$ is ball expanding on $U$. By normality, then, there is an open set $V$ such that $\Lambda\subseteq V\subset\cl{V}\subseteq U$.
\end{rem}

From \cite{Coven}, we know that tent maps have shadowing when the orbit of the critical point obeys certain rules. There are examples of tent maps where these rules are broken, yet the map still has shadowing on a subset of the interval. We make this idea precise in the following theorem, and follow it with an example of when this occurs.

\begin{thm}\label{thm:linear_shad}
Let $f\colon [0,1]\rightarrow [0,1]$ be a continuous function and let $\mathcal C$ be a closed, nowhere dense subset of $[0,1]$ that contains the critical points of $f$.  Suppose that $f$ is $C^1$ with gradient modulus strictly greater than 1 on every interval in $(0,1)\setminus \mathcal{C}$. If $\Lambda$ is a closed subset of $[0,1]$ that is disjoint from ${\mathcal C}\cup\{0,1\}$, then $f$ has shadowing, s-limit shadowing and h-shadowing on $\Lambda$.
\end{thm}
\ifproof\begin{proof}
Since $\Lambda$ and ${\mathcal C}\cup\{0,1\}$ are both closed, there is an open set $U\supseteq\Lambda$ with $\cl{U}\cap(\mathcal C \cup\{0,1\})=\emptyset$.  Notice that $U$ is a subset of only finitely many intervals, $U_1,\cdots, U_m$, of the (possibly countably many) intervals comprising $(0,1)\setminus\mathcal C$. To see this, note that otherwise we would be able to find an increasing or decreasing sequence of points $(x_n)$, where $x_{2k}\in U$ and $x_{2k+1}\in \mathcal C$, which would imply that that $\overline{U}$ and $\mathcal C$ have a common limit point. Clearly $f$ is expanding on $U'=U_1\cup\dots\cup U_m$ (take $\delta>0$ in the definition of expanding to be half the minimum of $\{d(x,c)\ :\ x\in\cl{U},\ c\in{\mathcal C}\}$, and $\mu>1$ to be the minimum gradient modulus of $f$ on the $U_i$, $i\le m$). By Remark \ref{expintmap}, $f$ is open on $U'$ and $f$ is ball expanding on $U'$, thus by Theorem \ref{CINEshadow}, $f$ has h-shadowing on $U'$. Shadowing follows directly from h-shadowing, and
Theorem \ref{h-shad_and_s-limit_shad} (\ref{h-shad_and_s-limit_shad1}) gives us that $f$ has s-limit shadowing on $\Lambda$, since $\Lambda\subseteq U'$.
\end{proof}\fi

\begin{rem}
If $f$ is piecewise linear map (i.e. there are is decomposition of $[0,1]$ into finitely many pieces, and $f$ is linear on each of them with gradient modulus strictly greater than 1),
then clearly assumptions of Theorem~\ref{thm:linear_shad} are satisfied. This extends to maps with infinitely many pieces of linearity as well.
\end{rem}

\begin{eg}\label{eg:tent_countereg} 
The tent map $T_2$ (with slope 2) is ball expanding (see Example \ref{eg:equiv_tent}) and, therefore, has h-shadowing by Theorem \ref{CINEshadow}.
As an example of a piecewise linear map which has shadowing on a subset of the interval $[0,1]$ but not on the interval itself, consider a tent map with gradient $\lambda\in(1,2)$ whose critical point $c$ is not recurrent. This map will not have shadowing \cite{Coven},  but by Theorem \ref{thm:linear_shad}, for any closed set $\Lambda\subseteq[0,1]$ for which $c\notin\Lambda$ we have that $f$ has shadowing (and s-limit shadowing and h-shadowing) on $\Lambda$.
\andysq
\end{eg}

The situation is less clear when the map in question is smooth and the conditions which imply shadowing in smooth maps are more subtle. We explore them in the next result.

\begin{eg}\label{eg:logistic_shad}
The logistic map $g_4(x)=4x(1-x)$ is conjugate to the tent map $T_2$, with gradient modulus 2, which we know to have shadowing \cite{Coven}, and as noted in Example \ref{eg:tent_countereg}, $T_2$ has h-shadowing. Thus we conclude that $g_4$ has h-shadowing (it is easy to show that h-shadowing is preserved between conjugate maps of a compact space). However neither Corollary \ref{col:uniexp_open_hshad} nor Theorem \ref{CINEshadow} apply here, since no member of the logistic family is expanding or ball expanding on any neighbourhood of the critical point $1/2$.
\andysq
\end{eg}

Recall that the Schwarzian derivative $S\big(f\big)(x)$ of a map $f$ at the point $x$ is given by
\[
S\big(f\big)(x)=\frac{f'''(x)}{f'(x)}-\frac{3}{2}\left(\frac{f''(x)}{f'(x)}\right)^2.
\]
 For $\eps>0$, a subset $N$ of a metric space $X$ is said to be an $\eps$-net if $N\cap B_\eps(y)\neq\nowt$ for all $y\in X$. We will interchangeably use the notation $D(f)=f'$, $D^{(2)}(f)=f''$, etc.

In \cite{misiurewicz}, Misiurewicz proves that for smooth maps with negative Schwarzian derivative and no sinks, there are open intervals $U$ on which the derivative of some iterate $m>0$ of the map is greater than 1, even where the derivative of the map itself is close to zero; this notion relies on the fact that $f^i(U)$ avoids neighbourhoods of the critical points for every $0\leq i<m$. Here we use a similar idea, proving that maps whose pre-critical points form a dense set have various shadowing properties on closed sets which do not contain their critical points.

\begin{thm} \label{misiurewicz type}
Suppose that $f\colon [0,1]\to[0,1]$ is continuous and $\Lambda\subset(0,1)$ is closed and strongly invariant. Suppose further that:
\begin{enumerate}
	\item \label{misiurewicz type: eps net} the set of critical points ${\mathcal C}$ is closed and for every $\eps>0$ there is some $m>0$ such that $f^{-m}({\mathcal C})$ is an $\eps$-net.
	\item $S\big(f(x)\big)$ is defined and non-positive for all $x\in [0,1]\setminus{\mathcal C}$;	
	\item $\Lambda\cap {\mathcal C}=\nowt$.
\end{enumerate}
Then $f$ has shadowing, s-limit shadowing and h-shadowing on $\Lambda$.
\end{thm}
\ifproof\begin{proof}
We claim that for some $m\in\N$, the derivative $D\big(f^m\big)(x)$ has absolute value strictly greater than 1, for all $x\in \Lambda$. Suppose that this is not the case, so that for every $m\in\N$, there is some $x_m\in \Lambda$ such that $\big|D\big(f^m\big)(x_m)\big|\le 1$. Let  $\d<d(\Lambda,{\mathcal C})/2$ and fix $m$ such that $f^{-m}({\mathcal C})$ forms a $\d/2$-net. Note that if $\d$ is sufficiently small then $x_m\in [\d/2,1-\d/2]$, because $x_m\in\Lambda$ which is a closed subset of $(0,1)$. We can find $c_1, c_2\in f^{-m}({\mathcal C})$ such that $c_1<x_m<c_2$. Since $\Lambda$ is invariant and disjoint from ${\mathcal C}$, $x_m\notin f^{-m}({\mathcal C})$.

Moreover as $f^{-m}(C)$ is closed, we may assume that $c_1=\max\{y\in f^{-m}({\mathcal C}):y<x_m\}$ and $c_2=\min\{y\in f^{-m}({\mathcal C}):x_m<y\}$. It follows that $c_2-c_1\le\d$ and $f^{-m}({\mathcal C})\cap (c_1,c_2)=\nowt$. Now the Schwarzian derivative $S(f^m)$ is well-defined on $(c_1,c_2)$ and by the proof of \cite[Proposition 11.3]{Devaney}, $S(f^m)\le0$ on $(c_1,c_2)$. It follows (see page 18 of \cite{misiurewicz}, for example) that$\big|D(f^m)\big|$ has no positive strict local minima on $(c_1,c_2)$.
Since $\big|D(f^m)(x_m)\big|\le 1$, we see that $\big|D(f^m)(y)\big|\le\big|D(f^m)(x_m)\big|\le 1$ for all $y$ in one of the intervals $(c_1,x_m)$ or $(x_m,c_2)$. Suppose, without loss of generality, that the first of these holds (the arguments for the second case are identical).  Then $|f^m(c_1)-f^m(x_m)|\le\big|D(f^m)(x_m)\big| \;|c_1-x_m|\le \d$, which contradicts the fact that $d(\Lambda,{\mathcal C})>\delta$. The proof of the claim is finished.

Since $\Lambda$ is compact, $\Lambda$ is disjoint from ${\mathcal C}$ and $S(f)$ exists on $[0,1]\setminus{\mathcal C}$, $Df$ is continuous on $[0,1]\setminus {\mathcal C}$,  it follows that there is some $\mu>1$ and a neighbourhood $U$ of $\Lambda$ on which the derivative of $f^m$ is greater than $\mu$ in absolute value. So, by Remark \ref{expintmap}, $f^m$ is open, expanding and ball expanding on a neighbourhood of $\Lambda$. So by Theorem \ref{CINEshadow} and Theorem \ref{iterateshadow} (which applies since $f(\Lambda)=\Lambda$), $f$ has h-shadowing on a neighbourhood of $\Lambda$; shadowing on $\Lambda$ is immediate, and s-limit shadowing follows from Theorem \ref{h-shad_and_s-limit_shad} (\ref{h-shad_and_s-limit_shad1}).
\end{proof}\fi

We end this section with an example which illustrates when Theorem \ref{misiurewicz type} can apply. In this example we draw heavily upon theory presented in \cite{ColletEckmann}. {While we provide precise references to facts used where possible, the reader is referred to \cite{ColletEckmann}
for many of the definitions, such as kneading sequence, $*$-product etc.

\begin{eg}\label{eg:tent_non_shad}
Consider the following sequence $K=RLLRRLRRRLRRRRL\ldots$ and observe that it is not recurrent under left shift map. Note that $K$ is not a $*$-product, i.e. $K\neq B*Q$ for any nonempty word $B$ and any sequence $Q\neq C$ (see \cite{ColletEckmann} for the full definition). Let $f_\mu\colon [-1,1]\rightarrow[-1,1]$ be the family of maps of the form $f_\mu(x)=1-\mu x^2$ with $\mu\in[1,2]$ and denote $J(f_\mu)=[f_\mu(1),1]$. Each $f_\mu$ is $C^3$, unimodal with critical point $0$, has maximum value $f_\mu(0)=1$, satisfies $f_\mu'(x)\neq 0$ and $S(f_\mu)(x)<0$ for every $x\neq 0$, and is surjective on $J(f_\mu)$.  Therefore each $f_\mu$ satisfies the definition of \emph{S-unimodal} from \cite{ColletEckmann}. Additionally $f_1(1)=0$, $f_2(1)=f_2(-1)=-1$ ($f_{\mu}$ is a so-called \textit{full family}) and so by \cite[III.1.2]{ColletEckmann} there is $\mu\in (1,2)$ such that $F=f_\mu$ has kneading sequence exactly equal to $K$ (note that $K$ is not the kneading sequence of either $f_1$ or $f_2$).

Since $K$ is infinite but not periodic, \cite[II.6.2]{ColletEckmann} implies that $F$ has no stable periodic orbit in $J(F)$. But in our case $F([-1,1])=J(F)$ so there is no stable periodic orbit of $F$ in $[-1,1]$ either. We can, therefore, apply \cite[II.5.5]{ColletEckmann} to see that the set of pre-critical points $\bigcup_{m\geq 0}F^{-m}(\set{0})$ is dense in $[-1,1]$. Since $K$ is not a $*$-product, by \cite[II.7.14 and II.7.12]{ColletEckmann} there is some $\lambda\in  (\sqrt{2}, 2)$ such that $F|_{J(F)}$ is conjugate to the map $g_\lambda|_{J(g_\lambda)}$ where $g_\lambda(x)=1- \lambda |x|$ for $x\in [-1,1]$. Repeating the arguments from the proof of \cite[Lemma~2]{Brucks} it is not hard to show that $g_\lambda$ restricted to $J(g_\lambda)$ is \emph{topologically exact}, that is for every open interval $I\subset J(g_\lambda)$ there is an $n>0$ such that $g_\lambda^n(I)\supset J(g_\lambda)$. By conjugacy, $F|_{J(F)}$ is also topologically exact, and since the pre-critical points are dense, for every interval $I\subset [-1,1]$ there is $n$ such that $F^m(I)\supset J(F)$ for every $m>n$. Now if we fix any $\eps>0$ and cover $[-1,1]$ with finitely many nondegenerate intervals $I_1,\ldots, I_k$ with diameters smaller than $\eps$ then there is also an $m>0$ such that $0\in J(F)\subset F^m(I_i)$ for each $i=1,\ldots,k$. This shows that condition \eqref{misiurewicz type: eps net} from Theorem~\ref{misiurewicz type} is satisfied.

Now we are going to show that $F$ does not have shadowing. Observe that the point $0$ is not recurrent, as in such a case $K$ has to be recurrent. Therefore we can find $\eps>0$ such that
\begin{align}
|F^n(0)-0|>\eps \mbox{ for every }n>0.\label{eg:tent_non_shad1}
\end{align}
and for which $F^2([-\eps,\eps])\subset[F(1),0)$. If $F$ had shadowing we would have some $\delta>0$ such that every $\delta$-pseudo-orbit is $\eps$-shadowed. We will show that this is not the case. Recall that $1<\mu<2$ and so $-1<F(1)<0$. Fix any $0<\delta<F(1)+1$ and let $x_0 = 0$, $x_1=1$, $x_2 = F(1)-\delta/2$ and for $i>2$, $x_i = F^{i-2}(x_2)$. Suppose that the $\delta$-pseudo-orbit $\Gamma = \{x_0,x_1,x_2,\ldots\}$ is $\eps$-shadowed by the point $y\in[-1,1]$. The pre-critical points are dense, so since $x_2<F^2(y)<0$ there is some least $n>2$ for which $0$ is between $F^n(y)$ and $x_n$. Notice that $F^2(0)\in[x_2,F^2(y)]$ but $0\notin[x_2,F^2(y)]$ so $F^3(0)$ lies between $F^3(y)$ and $F(x_2)=x_3$. Repeating this argument for every $i<n$ we get that $F^i(0)$ lies between $F^i(y)$ and $x_i$ but $0$ does not (by the definition of $n$). Thus $0$ and $F^n(0)$ both lie between $F^n(y)$ and $x_n$, and this contradicts \eqref{eg:tent_non_shad1} since $|F^n(y)-x_n|<\eps$. To finish our example, fix any closed and invariant set $\Lambda\subset (-1,0)\cup (0,1)$ and observe that by Theorem \ref{misiurewicz type} our map $F$ has shadowing, h-shadowing and s-limit shadowing on $\Lambda$, while it does not have (global) shadowing as shown above.

In other words, whilst not having any form of global shadowing, smooth maps such as $F$ have various shadowing properties on closed, invariant subsets of the interval which do not contain any critical point.

By changing the lengths of the blocks of $R$, one can easily show that there are uncountably many such interval maps.
\andysq
\end{eg}
}

\section{A Final Remark on $h$-shadowing Homeomorphisms}

Homeomorphisms with h-shadowing have some interesting properties, some of which we explore in this section.

Recall that a point $x\in X$ is an \emph{equicontinuity point} of $f$ if for every $\eps>0$ there is a $\delta>0$ such that given a point $y\in X$, $d(f^n(x),f^n(y))<\eps$ holds for all $n$ whenever $d(x,y)<\delta$. If every $x\in X$ is a point of equicontinuity then we say that $f$ is \emph{equicontinuous}.

\begin{thm}\label{eqcond:homeo}
Let $(X,d)$ be a compact metric space, and $f\colon X\rightarrow X$ be a homeomorphism. The following conditions are equivalent:
\begin{enumerate}
\item $f$ has h-shadowing;\label{eqcond:homeo:c1}
\item $f$ has shadowing and $f^{-1}$ is equicontinuous;\label{eqcond:homeo:c2}
\item $f$ is equicontinuous and $X$ is totally disconnected.\label{eqcond:homeo:c3}
\end{enumerate}
\end{thm}
\ifproof\begin{proof}
First we prove $\eqref{eqcond:homeo:c1}\Longrightarrow\eqref{eqcond:homeo:c2}$. Fix any $\eps>0$ and let $\delta>0$ be provided by h-shadowing. We may assume that $\delta<\eps$. Fix any $x,y\in X$ and assume that $d(x,y)<\delta$. For any $n>0$ the sequence
$$
f^{-n}(x),f^{-n+1}(x),\ldots, f^{-1}(x),y
$$
is a $\delta$-pseudo-orbit, so by h-shadowing of $f$ there exists $z$ such that
\[d(f^{-n+i}(x),f^{i}(z))<\eps
\]
for $i=0,\ldots, n$ and additionally $f^n(z)=y$. In other words
$$
d((f^{-1})^n(x),(f^{-1})^n(y))=d(f^{-n}(x),z)<\eps
$$
which proves that $f^{-1}$ is equicontinuous. Every map with h-shadowing has shadowing and so implication follows.

To prove $\eqref{eqcond:homeo:c2}\Longrightarrow\eqref{eqcond:homeo:c1}$, fix $\eps>0$ and choose $0<\gamma<\eps/2$ which satisfies the definition of equicontinuity for $\eps/2$. Let $0<\delta<\gamma$ be such that every $\delta$-pseudo-orbit is $\gamma$-shadowed. Fix any $\delta$-pseudo-orbit $x_0,x_1,\ldots, x_n$ and let $z$ be a point which $\gamma$-shadows it. By eqicontinuity of $f^{-1}$ we get that $d(f^{-i}(x_n),f^{-i}(f^n(z)))<\eps/2$ for every $i\geq 0$, since $d(x_n,f^n(z))<\gamma$.
Denote $y=f^{-n}(x_n)$
But then
\begin{eqnarray*}
d(f^j(y),x_j)&\leq &d(f^j(y),f^j(z))+d(f^j(z),x_j)\\
&\leq& d(f^{j-n}(x_n),f^{j}(z))+\gamma\\
&<&\eps/2+\eps/2
\end{eqnarray*}
and additionally $f^n(y)=x_n$.

To prove that $\eqref{eqcond:homeo:c2}\Longrightarrow\eqref{eqcond:homeo:c3}$ we argue as follows.

We claim that if $f^{-1}$ is equicontinuous, then every $x$ is recurrent under $f$, i.e. $x\in \w(x,f)$. To see this suppose that $d(x,\w(x,f))=\eps>0$. If $y\in \w(x,f)$, then, by the equicontinuity of $f^{-1}$, there is $\delta>0$ such that $d(f^{-i}(y),f^{-i}(z))<\eps/2$, for all $i$, provided $d(y,z)<\delta$. But there is $m>0$ such that $d(f^m(x),y)<\delta$, which is a contradiction.

We next claim that if $f^{-1}$ is equicontinuous and $f$ has shadowing, then $f$ is equicontinuous. Suppose that $x$ is not a point of equicontinuity for $f$ and pick $\eps>0$ such that for all $\delta>0$ there are $y$ and $m$ such that $d(x,y)<\delta$, but $d(f^m(y), f^m(x))\ge\eps$. For this $x$ and $\eps/4$, since $f^{-1}$ is equicontinuous we can choose $0<\eta<\eps/2$ such that $d(f^{-i}(x), f^{-i}(y))<\eps/4$, for all $i$, whenever $d(x,y)<\eta$. By shadowing of $f$, choose $0<\xi<\eta/3$ so that every finite $\xi$-pseudo orbit is $\eta/3$ traced. Now fix $y$ such that $d(x,y)<\xi$ and there is $m>0$ such that $d(f^m(x), f^m(y))\ge\eps$. Since $x$ and $y$ are recurrent there are $r>0$ and $s>0$ such that $d(x,f^r(x))<\xi$ and $d(y,f^s(y))<\xi$. Consider the two $\xi$-pseudo orbits of length $rs+1$ obtained by periodic concatenation of the finite orbits $\{x, f(x), \ldots , f^{r-1}(x)\}$ $s$ times followed by $x$ and  $\{y, f(y), \ldots , f^{s-1}(y)\}$ $r$ times followed by $y$. Let the orbit of $p$ $\eta/3$ trace the first of these and let $q$ $\eta/3$ trace the second. Note that $d(f^{sr+1}(p),x)<\eta/3<\eta$ and $d(f^{rs+1}(q),x)<\xi+\eta/3<\eta$. By the equicontinuity of $f^{-1}$, we have $d(f^m(p), f^{m-sr-1}(x))<\eps/4$ and $d(f^m(q), f^{m-sr-1}(x))<\eps/4$, so that
\begin{align*}
 d(f^m(x), f^m(y))&\le d(f^m(x), f^m(p))+d(f^m(p), f^{m-rs-1}(x))\\
 &\qquad+d(f^{m-rs-1}(x), f^m(q))+d(f^m(q),f^m(y))\\
&<\eta/3+\eps/4+\eps/4+\eta/3<\eps,
\end{align*}
which is a contradiction since $d(f^m(x), f^m(y))\ge\eps$.

Finally we show that $X$ is totally disconnected. Suppose not and let $x$ and $y$ be distinct points in a non-trivial connected component of $X$ with $d(x,y)=\eps>0$. Since $f$ is equicontinuous, there is $0<\eta<\eps/4$ such that $d(f^i(x),f^i(z))<\eps/4$, whenever $d(x,z)<\eta$. Since $f$ has shadowing, there is $\xi>0$ such that every $\xi$-pseudo orbit is $\eta/2$ traced. Since $x$ and $y$ are in the same connected component, there is a sequence of open $\xi/4$ balls, $B_0,\ldots, B_n$ such that $x\in B_0$ and $y\in B_n$ and $B_i\cap B_{i+1}\neq\emptyset$, for all $i$. Since every point in $X$ is recurrent, we can find $x_i\in B_i$ and $s_i>0$ such that $x_0=x$ and $x_m=y$ and $f^{s_i}(x_i)\in B_i$. Let $A_i$ denote the sequence $\{x_i, f(x_i), \ldots, f^{s_i-1}(x_i)\}$. Observe that the sequence $A_1A_1A_2A_2\dots A_nA_n$ followed by $y$ and the sequence $A_1A_2\dots A_nA_nA_{n-1}\dots A_1$ followed by $x$ are both $\xi$-pseudo orbits of the same length, $k$ say. Let the first of these be $\eta$ traced by the orbit of $p$ and the second be $\eta$ traced by the orbit of $q$. Since $d(x,p)<\eta/2$ and $d(x,q)<\eta/2$, we have
\begin{align*}
  d(x,y)&\le d(x,f^k(p))+d(f^k(p), f^k(x))+d(f^k(x),f^k(q))+d(f^k(q),y)\\
  &< \eta/2+\eps/4+\eps/4+\eta/2<\eps,
\end{align*}
which contradicts the fact that $d(x,y)\ge\eps$. This completes the proof of $\eqref{eqcond:homeo:c2}\Longrightarrow\eqref{eqcond:homeo:c3}$

Finally we prove that $\eqref{eqcond:homeo:c3}\Longrightarrow\eqref{eqcond:homeo:c2}$. Every equicontinuous map on a totally disconnected space has shadowing \cite[Prop. 4.7]{KurkaBook}. Since $f$ is a homeomorphism, every finite $\delta$-pseudo orbit for $f$ is a  $\delta'$-pseudo orbit for $f^{-1}$, where $\delta$ depends only on $\delta'$, so it is not hard to verify that $f$ has shadowing if and only if $f^{-1}$ has shadowing. Then by $\eqref{eqcond:homeo:c2}\Longrightarrow\eqref{eqcond:homeo:c3}$ $f^{-1}$ is equicontinuous, since $f$ is equicontinuous and $f^{-1}$ has shadowing.
\end{proof}\fi

\begin{rem}
Adding machines are a particular example of systems satisfying the conditions of Theorem \ref{eqcond:homeo}. Note that in \eqref{eqcond:homeo:c3} we cannot do better than totally disconnected and prove that $X$ is the Cantor set: the identity map on the union of Cantor set with some number of isolated points is equicontinuous and has shadowing.
\end{rem}

Recall that a map $f$ is \emph{(topologically) transitive} if for every pair of open sets $U$ and $V$, there is some $n\in\nat$ such that $f^n(U)\cap V\neq\emptyset$ and that $f$ is \emph{topologically mixing} if for every pair of open sets $U$ and $V$, there is an $N\in\nat$ for which $f^n(U)\cap V\neq\emptyset$ for every $n\geq N$. A map is said to be \emph{(topologically) weakly mixing} if the map $(f \times f)\colon (X \times X)\rightarrow (X \times X)$, defined by $(f \times f)(x,y)=(f(x),f(y))$, is transitive.

\begin{rem}\label{weakmixthm}
{Let $(X,d)$ be a compact metric space, and $f\colon X\rightarrow X$ be a topologically weakly mixing homeomorphism. It is easily seen (and known) that the maximal equicontinuous factor of a weakly mixing system is trivial, therefore if a weakly mixing homeomorphism has h-shadowing, X is a singleton. In other words, if $X$ has more than one element, then $f$ is not equicontinuous and by Theorem \ref{eqcond:homeo} $f$ does not have h-shadowing.}
\end{rem}

\begin{eg}\label{weakmixeg}
Consider any bi-infinite shifts of finite type with at least two elements. It has shadowing by \cite{Walters}, which demonstrates that h-shadowing and shadowing are not equivalent.
\andysq
\end{eg}


%
%
%
%
%

\bibliographystyle{plain}
\bibliography{BibtexW-limitsets}

\def\ocirc#1{\ifmmode\setbox0=\hbox{$#1$}\dimen0=\ht0 \advance\dimen0
  by1pt\rlap{\hbox to\wd0{\hss\raise\dimen0
  \hbox{\hskip.2em$\scriptscriptstyle\circ$}\hss}}#1\else {\accent"17 #1}\fi}
  \def\cprime{$'$}
\begin{thebibliography}{10}

\bibitem{AH}
N.~Aoki and K.~Hiraide.
\newblock {\em Topological theory of dynamical systems}, volume~52 of {\em
  North-Holland Mathematical Library}.
\newblock North-Holland Publishing Co., Amsterdam, 1994.
\newblock Recent advances.

\bibitem{BDG}
A.~D. Barwell, G.~Davies, and C~Good.
\newblock On the {$\omega$}-limit sets of tent maps.
\newblock {\em arXiv:1110.3219}, 2011.

\bibitem{Barwell}
A.~D. Barwell, C.~Good, R.~Knight, and B.~E. Raines.
\newblock A characterization of {$\omega$}-limit sets of shifts of finite type.
\newblock {\em Ergodic Theory and Dynamical Systems}, 30:21--31, 2010.

\bibitem{BGOR}
A.~D. Barwell, C.~Good, P.~Oprocha, and B.E. Raines.
\newblock Characterizations of {$\omega$}-limit sets in topologically
  hyperbolic systems.
\newblock {\em Discrete and Continuous Dynamical Systems}, 2012.

\bibitem{Blokh}
A.~Blokh, A.~M. Bruckner, P.~D. Humke, and J.~Sm{\'{\i}}tal.
\newblock The space of {$\omega$}-limit sets of a continuous map of the
  interval.
\newblock {\em Trans. Amer. Math. Soc.}, 348(4):1357--1372, 1996.

\bibitem{Bowen}
R.~Bowen.
\newblock {$\omega $}-limit sets for axiom {${\rm A}$} diffeomorphisms.
\newblock {\em J. Differential Equations}, 18(2):333--339, 1975.

\bibitem{Brucks}
K.~M. Brucks, B.~Diamond, M.~V. Otero-Espinar, and C.~Tresser.
\newblock Dense orbits of critical points for the tent map.
\newblock In {\em Continuum theory and dynamical systems ({A}rcata, {CA},
  1989)}, volume 117 of {\em Contemp. Math.}, pages 57--61. Amer. Math. Soc.,
  Providence, RI, 1991.

\bibitem{Chen}
L.~Chen.
\newblock Linking and the shadowing property for piecewise monotone maps.
\newblock {\em Proc. Amer. Math. Soc.}, 113(1):251--263, 1991.

\bibitem{ColletEckmann}
P.~Collet and J.-P. Eckmann.
\newblock {\em Iterated maps on the interval as dynamical systems}, volume~1 of
  {\em Progress in Physics}.
\newblock Birkh\"auser Boston, Mass., 1980.

\bibitem{Corless}
R.~M. Corless and S.~Y. Pilyugin.
\newblock Approximate and real trajectories for generic dynamical systems.
\newblock {\em J. Math. Anal. Appl.}, 189(2):409--423, 1995.

\bibitem{Corless2}
Robert~M. Corless.
\newblock Defect-controlled numerical methods and shadowing for chaotic
  differential equations.
\newblock {\em Phys. D}, 60(1-4):323--334, 1992.
\newblock Experimental mathematics: computational issues in nonlinear science
  (Los Alamos, NM, 1991).

\bibitem{Coven}
E.~M. Coven, I.~Kan, and J.~A. Yorke.
\newblock Pseudo-orbit shadowing in the family of tent maps.
\newblock {\em Trans. Amer. Math. Soc.}, 308(1):227--241, 1988.

\bibitem{Devaney}
R.~L. Devaney.
\newblock {\em An introduction to chaotic dynamical systems}.
\newblock Studies in Nonlinearity. Westview Press, Boulder, CO, 2003.
\newblock Reprint of the second (1989) edition.

\bibitem{FurBook}
H.~Furstenberg.
\newblock {\em Recurrence in ergodic theory and combinatorial number theory}.
\newblock Princeton University Press, Princeton, N.J., 1981.
\newblock M. B. Porter Lectures.

\bibitem{KO2}
M.~Kulczycki and P.~Oprocha.
\newblock Properties of dynamical systems with the aasp.
\newblock preprint.

\bibitem{KurkaBook}
Petr K{\ocirc{u}}rka.
\newblock {\em Topological and symbolic dynamics}, volume~11 of {\em Cours
  Sp\'ecialis\'es [Specialized Courses]}.
\newblock Soci\'et\'e Math\'ematique de France, Paris, 2003.

\bibitem{sakai}
K.~Lee and K.~Sakai.
\newblock Various shadowing properties and their equivalence.
\newblock {\em Discrete Contin. Dyn. Syst.}, 13(2):533--540, 2005.

\bibitem{misiurewicz}
M.~Misiurewicz.
\newblock Absolutely continuous measures for certain maps of an interval.
\newblock {\em Inst. Hautes \'Etudes Sci. Publ. Math.}, (53):17--51, 1981.

\bibitem{Nusse}
Helena~E. Nusse and James~A. Yorke.
\newblock Is every approximate trajectory of some process near an exact
  trajectory of a nearby process?
\newblock {\em Comm. Math. Phys.}, 114(3):363--379, 1988.

\bibitem{Pearson}
D.~W. Pearson.
\newblock Shadowing and prediction of dynamical systems.
\newblock {\em Math. Comput. Modelling}, 34(7-8):813--820, 2001.

\bibitem{Pennings}
Timothy Pennings and Jeffrey Van~Eeuwen.
\newblock Pseudo-orbit shadowing on the unit interval.
\newblock {\em Real Anal. Exchange}, 16(1):238--244, 1990/91.

\bibitem{Pil}
S.~Y. Pilyugin.
\newblock {\em Shadowing in dynamical systems}, volume 1706 of {\em Lecture
  Notes in Mathematics}.
\newblock Springer-Verlag, Berlin, 1999.

\bibitem{Urbanski}
F.~Przytycki and M.~Urba{\'n}ski.
\newblock {\em Conformal fractals: ergodic theory methods}, volume 371 of {\em
  London Mathematical Society Lecture Note Series}.
\newblock Cambridge University Press, Cambridge, 2010.

\bibitem{PEsakai}
K.~Sakai.
\newblock Various shadowing properties for positively expansive maps.
\newblock {\em Topology Appl.}, 131(1):15--31, 2003.

\bibitem{Walters}
Peter Walters.
\newblock On the pseudo-orbit tracing property and its relationship to
  stability.
\newblock In {\em The structure of attractors in dynamical systems ({P}roc.
  {C}onf., {N}orth {D}akota {S}tate {U}niv., {F}argo, {N}.{D}., 1977)}, volume
  668 of {\em Lecture Notes in Math.}, pages 231--244. Springer, Berlin, 1978.

\end{thebibliography}

\end{document}